\ifx\documentclass\undefined
\documentstyle[12pt]{article}
\else
\documentclass[12pt]{article}
\usepackage{latexsym}
\usepackage{amsfonts}
\usepackage{amsmath}
\usepackage{amssymb}
\fi

\author{K. Bezdek \thanks{Partially supported by the Hung.
Nat. Sci. Found (OTKA), grant no. NK 67867.} 
\thanks{Partially supported by a Natural Sciences and 
Engineering Research Council of Canada Discovery Grant.}
\and A. E. Litvak$\mbox{}^{\ddagger}$ }

\date{}

\font\tenBbb=msbm10 at 12pt         \font\sevenBbb=msbm9    \font\fiveBbb=msbm7
\newfam\Bbbfam
\textfont\Bbbfam=\tenBbb \scriptfont\Bbbfam=\sevenBbb
\scriptscriptfont\Bbbfam=\fiveBbb

\def\R{{\mathbb R}}

\def\kkk{\null\hfill $\Box$\smallskip}
\def\r{ \right}

\def\eps{\varepsilon}
\def\alp{\alpha}

\def\lam{\lambda}

\def\la{\left\langle}
\def\ra{\r\rangle}

\def\KK{{\bf K}}
\def\LL{{\bf L}}

\def\BB{{\bf B}}

\newcommand{\proof}{{\noindent\bf Proof:{\ \ }}}

\newtheorem{theorem}{Theorem}[section]

\newcommand{\vol}{{\rm vol}}

\newcommand{\Flt}{{\rm Flt}}
\newcommand{\wid}{{\rm w}}
\newcommand{\crv}{{\rm crv}}

\title{Covering convex bodies by cylinders and lattice points by flats 
}

\begin{document}

\maketitle

\begin{abstract}
In connection with an unsolved problem of Bang (1951) we give a lower bound for the sum 
of the base volumes of cylinders covering a $d$-dimensional convex body in terms of the 
relevant basic measures of the given convex body. As an application we establish lower 
bounds on the number of $k$-dimensional flats (i.e. translates of $k$-dimensional linear subspaces) needed to cover all the integer points of a given convex body in $d$-dimensional Euclidean space for $1\le k\le d-1$.
\end{abstract}

\section{Introduction}
\label{zero}

In a remarkable paper \cite{Ba} Bang has given an elegant proof of the plank conjecture of 
Tarski showing that if a convex body is covered by finitely many planks in $d$-dimensional 
Euclidean space, then the sum of the widths of the planks is at least as large as the minimal 
width of the body. A celebrated extension of Bang's theorem to $d$-dimensional normed spaces 
has been given by Ball in \cite{Ball2}. In his paper Bang raises also the important related 
question whether the sum of the base areas of finitely many cylinders covering a $3$-dimensional 
convex body is at least half of the minimum area of a 2-dimensional projection of the body. If true, 
then Bang's estimate is sharp due to a covering of a regular tetrahedron by two cylinders described 
in \cite{Ba}. We investigate this challenging problem of Bang in $d$-dimensional Euclidean space. 
Our main result is Theorem~\ref{covcyl} presented and proved in Section~\ref{two}. As a special case, we get that the sum of the base areas of finitely many cylinders covering a $3$-dimensional convex body is always at least one third of the minimum area $2$-dimensional projection of the body.

In \cite{BH} Bezdek and Hausel has established a discrete version of Tarski's plank problem by 
asking for the minimum number of hyperplanes that can cover the integer points within a convex 
body in $d$-dimensional Euclidean space. Theorem~\ref{Tal} of Section~\ref{four} gives an 
improvement of their result, which under some conditions improves also the corresponding 
estimate of Talata \cite{Ta}. A related but different problem of covering the lattice points 
within a convex body by linear subspaces was investigated in \cite{BHPT}. Last but not 
least, Theorem~\ref{covcyl} combined with some additional ideas leads to a lower bound on the number of $k$-dimensional flats (i.e. translates of $k$-dimensional linear subspaces) needed to cover all the integer points of a given convex 
body in $d$-dimensional Euclidean space for $1\le k\le d-1$. This is the topic of Section~\ref{tri} and its main 
result, Theorem~\ref{covlat}, actually improves the corresponding estimate of Talata \cite{Ta}.

\section{Notation}
\label{one}

In this paper we identify a $d$-dimensional affine space with $\R ^d$. 
By $|\cdot|$ and $\la \cdot , \cdot \ra$ we denote the canonical 
Euclidean norm and the canonical inner product on $\R ^d$. 
The canonical Euclidean ball and sphere in $\R^d$ are denoted by 
$\BB_2^d$ and $S^{d-1}$. By a subspace we always mean a linear subspace. 

By a convex body in $\R^d$ we always mean a compact convex set 
 with non-empty interior. The interior of $\KK$ is denoted by 
$\mbox{int} \KK$. Let $\KK\subset \R^d$ be a convex body
with the origin $0$ in its interior. We denote by 
$\KK^{\circ}$ the polar of $\KK$, i.e.
 $$
    \KK^{\circ} = \left\{ x \, \, | \, \, \la x, y \ra \leq 1 \, \,
    \mbox{ for every } \, \, y \in \KK \r\}.
 $$
The Minkowski functional of $\KK$ (or the gauge of $\KK$) is  
$$ 
  {\| x\|}_{\KK}=\inf \{\lambda >0 \ |\ x\in\lambda\KK\}. 
$$ 
If $\KK$ is a centrally symmetric convex body with its center of symmetry at the origin, then $\| x\|_{\KK}$
 defines a norm on $\R^d$ with the unit ball $\KK$. 

The {\it Banach-Mazur distance} between two convex  
bodies $\KK$ and $\LL$ in $\R ^d$ is defined by
$$  
    d(\KK, \LL) = \inf{\left\{ \lam > 0       \ \mid \ a\in \LL, \ b\in \KK, \ 
   \LL - a \subset T \left(\KK - b \r) \subset \lam \left(\LL - a\r) \r\}},
$$
where the infimum is taken over all linear operators $T : \R^d \to \R^d$. 
The Banach-Mazur distance between $\KK$ and the closed Euclidean 
ball $\BB_2^d$ (say, of unit radius) we denote by $d_{\KK}$. As it is well-known, 
John's Theorem (\cite{J}) implies that for every 
$\KK$, $d_{\KK}$ is bounded by $d$, while for centrally-symmetric convex body 
$\KK$, $d_{\KK}\leq \sqrt{d}$ (see e.g. \cite{Ball}). 

Given a convex body $\KK$ in $\R^d$ we denote its distance to symmetric 
bodies by 
\begin{equation} \label{sd}
   sd_{\KK} := \inf \left\{ \lam > 0 \ \mid \ a\in \R^d, 
   \ -(\KK-a) \subset \lam (\KK-a) \r\}. 
\end{equation}
Clearly, $sd_{\KK}\leq d_{\KK}\leq d$. In fact, $sd_{\KK}$ is one of the ways 
to measure the asymmetry of the convex body $\KK$. We refer to \cite{Gr} for 
the related discussion. 

Let $\KK$ be a convex body in $\R^d$. We denote its volume by $\vol (\KK)$. 
When we would like to emphasize that we take $d$-dimensional volume of a 
body in $\R^d$ we write $\vol _d (\KK)$.

Given a linear subspace (in short, a subspace) $E\subset \R^d$ we denote the orthogonal projection 
on $E$ by $P_E$ and the orthogonal complement of $E$ by $E^{\perp}$. 
We will use the following theorem, proved 
by Rogers and Shephard (\cite{RS}, see also \cite{Ch} and Lemma 8.8 in 
\cite{Pisier}).

\begin{theorem}\label{rsh} 
Let $1\leq k \leq d-1$. Let $\KK$ be a convex body in $\R^d$ and 
$E$ be a $k$-dimensional subspace of $\R^d$. Then   
$$
  \max _{x\in \R^d}\ \vol _{d-k} \left( \KK\cap \left( x + E^{\perp} 
  \r)\r) \vol _k (P_E\KK) \leq {d \choose k}   \vol _d (\KK). 
$$ 
\end{theorem}

\medskip

\noindent
{\bf Remark.} Note that the reverse estimate 
$$
  \max _{x\in \R^d}\ \vol _{d-k} \left( \KK\cap \left( x + E^{\perp} 
  \r)\r) \vol _k (P_E\KK) \geq   \vol _d (\KK) 
$$
is  a simple application of the Fubini Theorem and is correct for 
any measurable set $\KK$ in $\R^d$.

\medskip

We will be using the following parameters of a convex body 
$\KK$ with $0$ in its interior
$$
   M(\KK) := \int _{S^{d-1}} \|x\|_{\KK} \ d\sigma(x), 
$$
where $\sigma$ denotes the normalized Lebesgue measure on $S^{d-1}$,  
$M^*(\KK) := M(\KK ^{\circ})$, and 
$$
   MM^*(\KK) := \inf M(T(\KK-a))M^*(T(\KK-a)), 
$$
where the infimum is taken over all invertible linear maps 
$T: \R^d \to \R^d$ and all $a$ in the interior of $\KK$. 
Note that $M^*(\KK)$ is the half of mean width of $\KK$. 
Below we need the following theorem.

\begin{theorem}\label{mmm} There exist absolute positive 
constants $C$ and $\alp$ such that for every $d\geq 1$ and 
every convex body  $\KK$ in $\R^d$ one has 
$$
  MM^*(\KK) \leq C d^{1/3} \ln^{\alp}(d+1).
$$
Moreover, if $\KK$ is centrally symmetric then 
$$
  MM^*(\KK) \leq C \ln (d+1).
$$
\end{theorem}

The second estimate in this theorem is a well-known fact from Asymptotic 
Theory of finite dimensional normed spaces (see, e.g., \cite{Pisier, To}). 
In fact, it is a combination of results by Lewis (\cite{L}), by Figiel and 
Tomczak-Jaegermann (\cite{FT}) with a deep theorem by Pisier on the 
so-called Rademacher projection (\cite{Pi2}).  
The result in the general case is due to Rudelson (\cite{Rud}). 
The both estimates of the theorem plays an essential role in the 
Asymptotic Theory. 

The lattice width of a convex body $\KK$ in $\R^d$ is defined as
$$
  \wid (\KK, \mathbb Z ^d ) = \min\left\{ \max_{x\in \KK} \la x, y \ra - 
  \min_{x\in \KK} \la x, y \ra  \ \mid \ y\in \mathbb Z ^d, 
  \ y\ne 0 \r\}.
$$
Note that, if the origin is in the interior of $\KK$, then 
$$
 \wid (\KK, \mathbb Z ^d ) = \min\left\{ \|y\|_{\KK^{\circ}} + 
  \| - y\|_{\KK^{\circ}} \ \mid \ y\in \mathbb Z ^d, 
  \ y\ne 0 \r\}.
$$

The flatness parameter of $\KK$ is defined as 
$$
   \Flt (\KK) = \sup \wid (T\KK, \mathbb Z ^d ) , 
$$
where the supremum is taken over all invertible affine maps 
$\R^d \to \R^d$ satisfying $T\KK \cap \mathbb Z ^d = \emptyset$. 
The following theorem  was proved in \cite{Ban} for the 
centrally symmetric case and the case of an ellipsoid, and in 
\cite{BLPS} for the general case. It improves the previous bound by 
Kannan and Lov\'asz (\cite{KL}), who showed $\Flt (\KK) \leq C d^2$.

\begin{theorem}\label{flat}
There exist absolute positive 
constants $C$ and $c$ such that for every $d\geq 1$ and 
every convex body  $\KK$ in $\R^d$ one has 
$$
   c d \leq \Flt (\KK) \leq C d MM^* (\KK). 
$$
Moreover, $\Flt (\KK)\leq d$ if $\KK$ is an ellipsoid. 
\end{theorem}

\section{Covering by cylinders}
\label{two}

In this section we introduce a volumetric parameter related to 
covering by cylinders and provide corresponding estimates.

By a cylinder in $\R^d$ we always mean a 1-codimensional cylinder, that is, 
a set $C\subset \R^d$ that can be presented as $C = \ell + B$, where $\ell$ 
is a line containing $0$ in $\R^d$ and $B$ is a measurable set in 
$E: =\ell ^{\perp}$. Let $\KK \subset \R^d$ be a convex body and $C\subset \R^d$ 
be a cylinder. The cross-sectional volume of $C$ with respect to $\KK$ we denote by 
$$
 \crv _{\KK} (C) := \frac{\vol _{d-1} (C\cap E)}{\vol _{d-1} (P_E 
 \KK)} =\frac{\vol _{d-1} (P_E C)}{\vol _{d-1} (P_E 
 \KK)} = \frac{\vol _{d-1} (B)}{\vol _{d-1} (P_E \KK)} .
$$
It is easy to see that for every $(d-1)$-dimensional subspace $H\subset \R^d$ 
not containing $\ell$ one has 
$$
 \crv _{\KK} (C) = \frac{\vol _{d-1} (C\cap H)}{\vol _{d-1} (P \KK)}, 
$$
where $P$ is the projection on $H$ with the kernel $\ell$. 
We would also like to notice that
 for every invertible affine map $T: \R^d \to \R^d$ 
one has $\crv _{\KK} (C) = \crv _{T \KK} (T C)$.

\begin{theorem}\label{covcyl}
Let $\KK$ be a convex body in $\R^d$. Let $C_1, \dots, C_N$ 
be cylinders in $\R^d$ such that 
$$
  \KK\subset \bigcup _{i=1}^N C_i .
$$
Then 
$$
    \sum _{i=1}^N  \mbox{\rm crv} _{\KK} (C_i) \geq \frac{1}{d} .
$$ 
Moreover, if $\KK$ is an ellipsoid then 
$$
    \sum _{i=1}^N  \mbox{\rm crv} _{\KK} (C_i) \geq 1.  
$$
\end{theorem}

\medskip

\proof In this proof we denote $v_{n} := \vol _{n} (\BB _2^{n})$. 
Every $C_i$ can be presented as $C_i = \ell _i + B_i$, 
where $\ell _i$ is a line containing $0$ in $\R^d$ 
and $B_i$ is a body in $E_i: =\ell _i^{\perp}$.

We first prove the theorem for ellipsoids. Since 
$\crv _{\KK} (C) = \crv _{T \KK} (T C)$ for every invertible affine map 
$T: \R^d \to \R^d$, we may assume that $\KK = \BB_2^d$.  Then 
$$
  \crv _{\KK} (C_i) = \frac{\vol _{d-1} (B_i)}{v_{d-1}}.
$$

Consider the following (density) function on $\R^d$
$$
  p(x) = 1/\sqrt{1-|x|^2}
$$
for $|x|<1$ and $p(x)=0$ otherwise. The corresponding 
measure on $\R^d$ we denote by $\mu$, that is $d\mu (x) = p(x) dx$.  
Let $\ell$ be a line  containing $0$ in $\R^d$ and $E=\ell ^{\perp}$. 
It follows from direct calculations that for every $z\in E$ with $|z| < 1$ 
$$
  \int _{\ell +z} p(x) \ dx =  \pi.
$$
Thus, we have 
$$
  \mu (\BB_2^d) = \int _{\BB_2^d} p(x) \ dx = \int _{\BB_2^d \cap E} 
  \int _{\ell +z} p(x) \ dx\ dz = \pi \ v_{d-1} 
$$
and for every $i\leq N$ 
$$
  \mu (C_i) = \int _{C_i} p(x)\ dx = \int _{B_i} \int _{\ell _i +z} 
  p(x) \ dx\ dz = \pi \ \vol_{d-1} \left(B_i\r).  
$$
Since $\BB_2^d \subset \bigcup _{i=1}^N C_i $, we obtain 
$$
  \pi \ v_{d-1} =\mu (\BB_2^d) \leq \mu\left(\bigcup _{i=1}^N C_i\r) 
  \leq \sum _{i=1}^N  \mu\left( C_i\r) =  \sum _{i=1}^N   
  \pi \ \vol_{d-1} \left(B_i\r).
$$
It implies
\begin{equation}\label{ellcyl}
  \sum _{i=1}^N \crv _{\BB_2^d} (C_i) = \sum _{i=1}^N  
  \frac{\vol _{d-1} (B_i)}{v_{d-1}} \geq 1 .
\end{equation}

\smallskip

Now, we show the general case. For $i\leq N$  
 denote  $\bar C_i = C_i \cap  \KK$ and note that 
$$
  \KK\subset \bigcup _{i=1}^N \bar C_i \quad \quad  \mbox{ and }
  \quad \quad 
  P_{E_i} \bar C_i = B_i \cap P_{E_i} \KK . 
$$
Since $\bar C_i \subset \KK$  we have also 
$$
  \max _{x\in \R^d}\  \vol _{1}  \left( \bar C_i \cap 
  \left( x + \ell _i \r) \r) \leq 
  \max _{x\in \R^d}\  \vol _{1}  \left( \KK \cap 
  \left( x + \ell _i \r) \r) .
$$
Therefore, applying Theorem~\ref{rsh} (and Remark after it, saying that 
we don't need convexity of $\bar C_i$)  we obtain for every $i\leq N$ 
$$
  \mbox{\rm crv} _{\KK} (C_i) = 
  \frac{\vol _{d-1} (B_i)}{\vol _{d-1} (P_{E_i} \KK)} \geq  
  \frac{\vol _{d-1} (P_{E_i} \bar C_i)}{\vol _{d-1} (P_{E_i} \KK)} 
$$
$$
  \geq  \frac{\vol _{d} (\bar C_i) }{\max _{x\in \R^d}\  \vol _{1}  \left( \bar C_i \cap 
  \left( x + \ell _i \r) \r)  } \ 
  \frac{\max _{x\in \R^d}\  \vol _{1}  \left(\KK \cap 
  \left( x + \ell _i \r) \r) }{d \vol _{d} (\KK)} \geq \frac{\vol _{d} (\bar C_i)  }{d \vol _{d} (\KK)}.  
$$
Using that $\bar C_i$'s covers $\KK$, we observe 
$$
 \sum _{i=1}^N  \mbox{\rm crv} _{\KK} (C_i) 
 \geq \frac{ 1 }{d},   
$$
which completes the proof. 
\kkk

\medskip

\noindent 
{\bf Remark 1. } If $\KK$ is close to the Euclidean ball 
(and $d$ is not very big), then the following estimate can 
be better than the general one 
$$
 \sum _{i=1}^N  \mbox{\rm crv} _{\KK} (C_i) 
 \geq \frac{ 1 }{d_{\KK}^{d-1}}.   
$$
It can be obtained as follows: 
Using that 
$\crv _{\KK} (C) = \crv _{T \KK} (T C)$ for an invertible affine 
transformation, we may assume that $\BB_2^d$ is a distance ellipsoid 
for $\KK$, namely assume that $\BB_2^d \subset \KK\subset d_{\KK} \BB_2^d$. 
Then 
$$
  \sum _{i=1}^N \crv _{\KK} (C_i) = \sum _{i=1}^N  \frac{\vol _{d-1} (B_i)}{
  \vol _{d-1} (P_{E_i} \KK)} \geq \sum _{i=1}^N \frac{\vol _{d-1} (B_i)}{
  \vol _{d-1} (P_{E_i} d_{\KK} \BB_2^d)} 
$$
$$
  \geq d_{\KK}^{-d+1} \ \sum _{i=1}^N \crv _{\BB_2^d} (C_i) \geq d_{\KK}^{-d+1}
$$
 (in the last inequality we used ``moreover" part of Theorem~\ref{covcyl}).
Recall  that $d_{\KK}\leq \sqrt{d}$ for any centrally 
symmetric convex body $\KK$ in $\R^d$ and $d_{\KK}\leq d$ in general. Thus, 
if $d=3$ and $\KK$ is a centrally-symmetric convex body close to the Euclidean ball, then this estimate is better than
the general one given by Theorem~\ref{covcyl}.

\medskip 
\noindent 
{\bf Remark 2. } Note that the proof of Theorem~\ref{covcyl} can be extended 
to the case of cylinders of other dimensions. Indeed, given $k<d$ 
define a $k$-codimensional cylinder $C$ as a set which can be presented in the 
form $C = H + B$, where $H$ is a $k$-dimensional subspace of $\R^d$ 
and  $B$ is a measurable set in $E: =H^{\perp}$.  As before, given  a convex body $\KK$ 
and a $k$-codimensional cylinder $C= H + B$ denote 
$$
 \crv _{\KK} (C) := \frac{\vol _{d-k} (C\cap E)}{\vol _{d-k} (P_E 
 \KK)} =\frac{\vol _{d-k} (P_E C)}{\vol _{d-k} (P_E 
 \KK)} = \frac{\vol _{d-k} (B)}{\vol _{d-k} (P_E \KK)} .
$$
Repeating the proof of Theorem~\ref{covcyl} (the general case), we obtain that if 
a convex body $\KK$ is covered by $k$-codimensional cylinders $C_1$, \ldots, $C_n$, 
then 
$$
    \sum _{i=1}^N  \mbox{\rm crv} _{\KK} (C_i) \geq \frac{1}{{d \choose k}} .
$$
As was noted by Bang (\cite{Ba}), the case $k=d-1$ here corresponds to the 
``plank problem", indeed, in this case we have the sum of relative widths of the body. 
As we mentioned in the introduction, Ball (\cite{Ball2}) proved that such sum should exceed 
1 in the case of centrally symmetric body $\KK$, while the general case is still open. 
Our estimate implies the lower bound $1/d$. Of course, Ball's Theorem implies   
the estimate $1/ sd_{\KK}$.

\section{Covering lattice points by lines and flats}
\label{tri}

\begin{theorem} \label{covlat}
Let $\KK$ be a convex body in $\R^d$ containing the origin in its 
interior. Let $\ell_1, \dots, \ell_N$ be lines in $\R^d$ such that 
$$
  \KK \cap \mathbb Z ^d \subset  \bigcup _{i=1}^N \ell_i .
$$
Then 
$$
    N \geq    \left( \frac{ \wid \left( \KK\cap -\KK, 
   \mathbb Z ^d \r) }{C d \ MM^* \left(\KK\cap -\KK \r) } \r) ^{d-1}
   \geq   \left( \frac{ \wid \left( \KK\cap -\KK, 
   \mathbb Z ^d \r) }{C_0 d \ \ln(d+1) } \r) ^{d-1}, 
$$
where $C$ and $C_0$ are absolute positive constants. If, in addition, 
 $-\KK\subset sd_{\KK} \KK$ (that is, if infimum in (\ref{sd}) 
attains at $a=0$), then 
$$
   N  \geq  \left( \frac{ \wid \left( \KK, \mathbb Z ^d \r) 
   }{ C \ sd_{\KK}\ d \ MM^* \left(\KK \r)} \r) ^{d-1} 
   \geq \left( \frac{ \wid \left( \KK, \mathbb Z ^d \r) 
   }{ C_0 \  d^{7/3} \ln^{\alp}(d+1) }  \r) ^{d-1},  
$$
where $C$, $C_0$, and $\alp$ are absolute positive constants. 

Moreover, if $\KK$ is an ellipsoid centered at the origin, then 
$$
   N \geq \left( \frac{ \wid \left( \KK, 
   \mathbb Z ^d \r) }{2 d  } \r) ^{d-1}. 
$$
\end{theorem}

\proof Let $\lam>0$ be such that 
$$
  \KK\subset \bigcup _{i=1}^N \left(\ell_i + \lam \KK\r) 
  \quad \quad \mbox{ and } \quad \quad 
  \KK\not\subset \bigcup _{i=1}^N \left(\ell_i + \lam \mbox{ int} \KK\r) .
$$
Since $0\in \KK$, we have $0\in l_i$ for some $i$, which clearly implies 
that $\lam \leq 1$.

For $i\leq N$ let $H_i$ denote the $(d-1)$-dimensional subspace orthogonal to 
$\ell _i$ and let $P_i$ denote the orthogonal projection on $H_i$. 
We define 
$$
   C_i :=\ell_i + \lam \KK = \ell_i + \lam P_i \KK.  
$$
Then $\crv _{\KK} (C_i) = \lam ^{d-1}$. 
Theorem~\ref{covcyl} implies $N\geq c^d \lam ^{-d+1}$, where 
$c$ is a positive absolute constant.

Now,  $\KK\not\subset \cup _{i=1}^N \left(\ell_i + \lam \mbox{ int}\KK\r)$ 
if and only if there exists $x\in \KK$ such that for every $i\leq N$ one has 
$x\not\in  \ell_i + \lam \mbox{ int}\KK$, i.e. 
$\left(x-\lam \mbox{ int}\KK\r) \cap \ell_i =  \emptyset$. Let $y=(1-\lam/2)x$. 
By convexity of $\KK$ we have 
$$
  \left( y+\frac{\lam}{2}\ \left(\KK\cap -\mbox{ int}\KK \r)\r) \subset \KK 
  \cap \left(x -\lam \mbox{ int}\KK \r). 
$$
Since $\KK \cap \mathbb Z ^d \subset  \cup _{i=1}^N \ell_i$, we obtain 
$$
 \left( y+\frac{\lam}{2}\ \left(\KK\cap -\mbox{ int}\KK \r) \r)\cap 
 \mathbb Z ^d =  \emptyset .
$$ 
 Using  Theorem~\ref{flat} (and, if needed, approximating $\lam$ by  
$\lam - \eps$ with small enough $\eps$), we observe 
$$
   \frac{\lam}{2}\ \wid \left(\KK\cap -\KK, \mathbb Z ^d \r) = 
   \wid \left(y+\frac{\lam}{2}\ \left(\KK\cap -\KK \r), \mathbb Z ^d \r) 
$$
$$
  \leq \Flt (\KK\cap -\KK )  \leq  C d \ MM^* \left(\KK\cap -\KK \r), 
$$
where $C$ is an absolute constant. 
Thus, 
$$
   N\geq c^d  \lam ^{-d+1} \geq c^d   
   \left( \frac{ \wid \left( \KK\cap -\KK, 
   \mathbb Z ^d \r) }{2 C d \ MM^* \left(\KK\cap -\KK \r) } \r) ^{d-1}.
$$
This shows the left-hand side of the first estimate.  
The right-hand side follows by Theorem~\ref{mmm}. Note that in the case of ellipsoid 
we have $C = c = 1$, $MM^* \left(\KK\cap -\KK \r)=1$, which implies 
the ``moreover" part of the theorem.  

The second estimate follows the same lines. For the sake of 
completeness we sketch it. Let $0<\lam \leq sd_{\KK}$ be such that 
$$
  \KK\subset \bigcup _{i=1}^N \left(\ell_i - 2\lam \KK\r) 
  \quad \quad \mbox{ and } \quad \quad 
  \KK\not\subset \bigcup _{i=1}^N \left(\ell_i - \lam \mbox{ int}\KK\r) .
$$
Repeating arguments of the first part we obtain that  
$N\geq c^d  \lam ^{-d+1}$ and 
$\left(x + \lam \mbox{ int}\KK\r) \cap \ell_i =  \emptyset$ for every $i\leq N$. 
Convexity of $\KK$ and the inclusion $-\KK \subset sd_{\KK} \KK$
yields for $y=(1-\lam/(sd_{\KK}+1))x$ 
$$
  \left( y+\frac{\lam}{sd_{\KK}+1}\ \mbox{ int}\KK \r) \subset \KK \cap 
  \left(x + \lam \mbox{ int}\KK \r). 
$$
It implies 
$$
 \left( y+\frac{\lam}{sd_{\KK}+1}\ \mbox{ int}\KK \r)\cap 
 \mathbb Z ^d =  \emptyset 
$$ 
 and, by Theorem~\ref{flat}, 
$$
   \frac{\lam}{sd_{\KK}+1}\ \wid \left(\KK, \mathbb Z ^d \r) \leq 
   C_1 d \ MM^* \left(\KK \r).
$$
Therefore,
$$
   N\geq c^d  \lam ^{-d+1} \geq  c^d  
   \left( \frac{ \wid \left( \KK, \mathbb Z ^d \r) }{C_1 \left( 
   sd_{\KK}+1\r) d \ MM^* \left(\KK \r)} \r) ^{d-1}, 
$$
which proves the left-hand estimate (with $C=2 C_1$). Since $sd_{\KK}\leq d$, 
Theorem~\ref{mmm} implies the right-hand side inequality. 
\kkk

\medskip

\noindent
{\bf Remark. } It is not difficult to see that the proof above can be extended almost 
without changes to the case of $k$-dimensional flats instead of lines (one needs 
to use Remark 2 following Theorem~\ref{covcyl}). In particular, for a centrally symmetric 
body $\KK=-\KK$, whose integer points are covered by the $k$-dimensional flats $H_1$, 
\ldots, $H_N$ we have 
$$
   N  \geq  \left( \frac{\wid \left( \KK, \mathbb Z ^d \r) (d-k)  
   }{ C \  d^2 \ \ln(d+1) } \r) ^{d-k} .
$$
We omit the details and precise estimates in the non-symmetric case. 

\section{Covering lattice points by hyperplanes}
\label{four}

The following theorem improves the estimate of the Remark after Theorem~\ref{covlat} 
in the case $k=d-1$, extending a Bezdek-Hausel result from \cite{BH}.

\begin{theorem}\label{Tal} Let $\KK$ be a centrally symmetric 
(with respect to the origin) convex body in $\R^d$. 
Let $H_1, \dots, H_N$ be hyperplanes in $\R^d$ such that  
$$
   \KK \cap  \mathbb Z ^d \subset \bigcup _{i=1}^N H_i. 
$$
Then 
$$
   N\geq c \ \frac{\wid (\KK, \mathbb Z ^d)}{d \ MM^*(\KK)}
   \geq c_0 \ \frac{\wid (\KK, \mathbb Z ^d)}{d \ \ln(d+1)},  
$$
where $c$, $c_0$ are absolute positive constants. 
\end{theorem}

\proof 
The proof is based on the Ball's solution of the plank problem. 
Namely, we use that given a centrally symmetric body $\KK\subset \R^d$ 
and $N$ hyperplanes $H_1, \ldots, H_N$ in $\R^d$ there exists $x\in R^d$ such 
that 
$$
  \LL: = x+ \frac{1}{N+1}\ \KK \subset \KK
$$ 
and the interior of $\LL$  is not met by any $H_i$
(see Corollary or abstract in \cite{Ball2}). 

Since all integer points of $\KK$ are covered by $H_i$'s, we observe that 
$$
  \mbox{ int} \LL \cap  \mathbb Z ^d = \emptyset.   
$$
 Applying Theorem~\ref{flat}, we obtain 
$$
 \frac{1}{N+1}\ \wid \left(\KK, \mathbb Z ^d \r) = 
 \wid \left(\LL, \mathbb Z ^d \r) \leq \Flt (\KK) \leq  
 C d \ MM^* \left(\KK \r), 
$$
where $C$ is an absolute constant. Together with Theorem~\ref{mmm} 
it implies the desired result. 
\kkk

\vspace{1cm}

\medskip

\noindent
K\'aroly Bezdek,
Department of Mathematics and Statistics,
2500 University drive N.W.,
University of Calgary, AB, Canada, T2N 1N4.
\newline
{\sf e-mail: bezdek@math.ucalgary.ca}
 
\smallskip

\noindent
A.E. Litvak,
Department  of Mathematical and Statistical Sciences,
University of Alberta, Edmonton, AB, Canada, T6G 2G1.
\newline
{\sf e-mail: alexandr@math.ualberta.ca}

\end{document}